\documentclass{tsp-short}

\usepackage{bbm}
\newtheorem{thm}{Theorem}

\theoremstyle{remark}
\newtheorem{rem}{Remark}
\theoremstyle{definition}
\newtheorem{dfn}{Definition}

\newcommand{\Prob}{\mathrm{P}}
\newcommand{\E}{\mathrm{E}}

\newcommand{\e}{\mathrm{e}}

\newcommand{\ve}{\varepsilon}
\newcommand{\vf}{\varphi}

\begin{document}

\title
{Iterated Logarithm Law for Sizes of Clusters in Arratia flow}

\author{A. A. Dorogovtsev}
\address{Institute of Mathematics of the NAS of Ukraine}
\email{adoro@imath.kiev.ua}

\author{A. V. Gnedin}
\address{Queen Mary, University of London}
\email{A.Gnedin@qmul.ac.uk}

\author{M. B. Vovchanskii}
\address{Institute of Mathematics of the NAS of Ukraine}
\email{vovchansky.m@gmail.com}

\subjclass[2000]{Primary 60J65; Secondary 60D05}
\keywords{Arratia flow, cluster, Brownian motion, Gaussian processes, concentration of measure}

\begin{abstract}
The asymptotics of sizes of clusters for the Arratia flow is considered, the Arratia flow being a system of coalescing Wiener processes starting from the real axis and independent before they meet. A cluster at time $t$ is defined as a set of particles that have glued together not later than at $t.$  The results obtained are remarked to hold for any Arratia flow with a Lipschitz drift.
\end{abstract}

\maketitle

In the paper, we study internal properties of  the Arratia flow. The Arratia flow has been introduced  in~\cite{arratia_phd} and is an instance of Harris flows~\cite{harris,warren}. Informally speaking, it is a system of coalescing Wiener processes such that each pair of them has coordinates that are independent before a collision appears and then merge together. As is known~\cite{harris}, all particles within the flow glue into a finite number of clusters for any nonzero moment of time. We ask about the behavior of such a cluster, namely, one that contains a particle starting from $0,$ when time is small.

The next definition may be found in~\cite{harris,warren}.
\begin{dfn}
A family of mappings $\{y(u,\cdot)\mid u\in\mathbb{R}\},$ with $y(u)\equiv \{y(u,t)| t\in \mathbb{R}^{+}\},$ is called the Arratia flow if

1) for any $u,$ \ \ $y(u)$ is a Wiener process; $y(u,0)=u;$

2) for $u_1\leq u_2,$ \ \ $y(u_1)\leq y(u_2);$

3) the joint characteristic of the martingale parts of $y(u_1)$ and $y(u_2)$ equals
$$
\int^t_0\mathbbm{1}_{\{y(u_1,s)=y(u_2,s)\}}ds.
$$
\end{dfn}

Consider $A_t=\{u\mid y(0,t)=y(u,t)\}.$ We say that $A_t$ is a cluster formed by all particles that have glued with the particle from $0$ till the moment $t.$ Despite the supremum of $A_t$ is questioned to belong to $A_t,$ the size of a cluster can be defined as
\begin{equation}
\label{eq:nu_definition}
 \nu(t)=\lambda \{u\mid y(0,t)=y(u,t)\},
\end{equation}
where $\lambda$ is the Lebesgue measure on the line. In fact, it is shown in~\cite{dorogovtsev:continuity} that any Arratia flow has a version such that a mapping $u\to \{y(u,t)\mid t\in [0;1]\}$ is right continuous in $C([0;1]).$ This version being chosen, it implies $A_t$ is well defined, at least for $t\le 1.$ For the sake of clearance, we suppose hereinafter that such version of the flow is chosen.

The distribution of $\nu(t)$ may be written down in terms of collision times, whose notation  will be used in the sequel. On this way, for $f,g\in \mathrm{C}(\mathbbm{R}^{+}),$ define 
collision time
$$
\tau[f,g]=\inf\{t\mid f(t)=g(t)\},
$$
the infimum over the empty set being defined as infinity. Also define
$$
\theta(y)=\inf\{s\mid W(s)=y\},
$$
where $W$ is a standard Wiener process started from $0.$ Then, with the monotonicity property of an Arratia flow, it directly implies that
\begin{multline}
 \label{eq:1}
\Prob\{\nu(t)\geq r\}= P\{\forall u\in(0; r]: \tau[y(0), y(u)]\leq t\}
=\Prob\{\tau[y(0), y(r)]\leq t\}= \\
  =\Prob\left\{\theta\left(\frac{r}{\sqrt{2}}\right) \leq t\right\},
\end{multline}
since the difference $\frac{y(u)-y(0)}{\sqrt{2}}$ is a Wiener process until the collision happens.

In~\cite{darling_sizez_clusters}, the asymptotics of sizes of  clusters as $t$ grows to $\infty$ were obtained  in the more general  case of Harris flows. More precisely, consider a Harris flow~\cite{harris} with an infinitesimal covariance function $\psi.$ Here, $\psi$ is supposed to be non-negative definite, continuous on the real line, and satisfying the Lipschitz condition outside each interval $(-c;c),\ c> 0,$ and its  spectral distribution is not of a pure jump type;
$\psi(0)=1.$ Define $\nu(t)$ as the size of the cluster containing a particle that starts from $0$ in the same way as it has been done for an Arratia flow. The following theorem is taken from~\cite{darling_sizez_clusters}.

\begin{thm}
 \label{theor:darling}
Consider a Harris flow $\{X(u,t)\mid u\in\mathbb{R}, t\in
\mathbb{R^+}\}.$ Then, for each $y>0,$
$
\lim_{t\rightarrow\infty}\Prob\{\frac{\nu(t)}{\sqrt{t}}<y\}$ exists and equals  $\Prob\{\inf_{s\in{[0;2]}} W_y(s)> 0\},$ where $W_y$ is a standard Wiener process starting from $y.$

\end{thm}

It shows that when time goes to infinity, the asymptotic behaviour of a cluster inside any Harris flow is the same as that in the case of an Arratia flow. Another set of results concerning the divergence of particles from the start points may be found in~\cite{shamov}.

The main result of the paper is the following.

\begin{thm}
\label{theor:sizes_clusters_arratia} Let $y$ be an Arratia flow, and let $\nu$ be defined via~\eqref{eq:nu_definition}. Then a.s. 
$$
{\limsup}_{t\rightarrow 0+} \frac{\nu(t)}{\sqrt{2t\ln\ln t^{-1}}}\geq 1,
$$
$$
{\limsup}_{t\rightarrow 0+} \frac{\nu(t)}{2\sqrt{t\ln\ln t^{-1}}}\leq 1.
$$
\end{thm}

\proof  First, we prove the estimate from above. Define a function
$\vf\colon t\mapsto 2\sqrt{t\ln\ln t^{-1}}.$  Fix $\ve >0.$ For $n\in\mathbbm{N},$ put $t_n=\alpha^{n}, \alpha\in (0;1),$ and define $A_n=\left\{\frac{\nu(t_n)}{\vf(t_n)}\ge 1+\ve\right\}.$ For any $\alpha$ starting from some $n_0,$ the variables $(1+\ve)\vf(t_n)$ are well defined.  We restrict ourselves in what follows to only $n$ such that $u_n=(1+\ve)\vf(t_n)$ exists. Then, by \eqref{eq:1},
\begin{equation*}
 \Prob(A_n)=\Prob\{\tau[y(u_n),y(0)]\leq t_n\}=\Prob\left\{\theta(\frac{u_n}{\sqrt{2}})\leq t_n\right\}.
\end{equation*}
Thus,
\begin{multline*}
 \sum_{n\geq 1}\Prob(A_n)=\sum_{n\geq 1} \sqrt{\frac{2}{\pi}} \int_{\frac{u_n}{\sqrt{2t_n}}}^{+\infty}  \e^{-\frac{v^2}{2}}dv \le \frac{2}{\sqrt{\pi}} \sum_{n\ge 1} \frac{\sqrt{t_n}}{u_n} \e^{-\frac{u_n^2}{4t_n}}=
\\
=\frac{1}{\sqrt{\pi}(1+\ve)} \sum_{n\ge 1} \frac{1}{\sqrt{\ln\ln t_n^{-1}}} \e^{-(1+\ve)^2 \ln\ln t_n^{-1}}
=\frac{1}{\sqrt{\pi}(1+\ve)} \cdot
\\
\cdot\sum_{n\ge 1} \frac{1}{\sqrt{\ln (n\ln\alpha^{-1})}} \e^{-{(1+\ve)^2} \ln (n\ln\alpha^{-1}) }<+\infty.
\end{multline*}
Therefore, by the Borel-Cantelli lemma,
\begin{equation}
\label{eq:2}
 \Prob\{\mathrm{starting\ from\ some\ number\ } n\quad \frac{\nu(t_n)}{\varphi(t_n)}<1+\ve \}=1.
\end{equation}
For $t\in [t_{n+1};t_n],$
\begin{equation}
\label{eq:3}
 \frac{\nu(t)}{\vf(t)}\leq \frac{\nu(t_{n})}{\vf(t_{n+1})}=\frac{\nu(t_{n})}{\vf(t_{n})}\cdot \frac{\vf(t_n)}{\vf(t_{n+1})}=\frac{\nu(t_{n})}{\vf(t_{n})}\cdot \sqrt{\frac{ \ln (n\ln\alpha^{-1})}{\alpha \ln ((n+1)\ln\alpha^{-1})}}<\frac{\nu(t_{n})}{\vf(t_{n})}\cdot \alpha^{-\frac{1}{2}}.
\end{equation}
Equations \eqref{eq:2} and \eqref{eq:3} yield
\begin{equation*}
 \Prob\{{\sup \lim}_{t\rightarrow 0+} \frac{\nu(t)}{2\sqrt{t\ln\ln t^{-1}}}\leq (1+\ve)\alpha^{-\frac{1}{2}} \}=1.
\end{equation*}
Since $\ve\geq 0,$ and $\alpha\in(0;1)$ are arbitrary, we obtain
$$
\Prob\{{\sup \lim}_{t\rightarrow 0+} \frac{\nu(t)}{2\sqrt{t\ln\ln t^{-1}}}\leq 1 \}=1.
$$

In order to prove the first estimate in the statement of the theorem, we need the following two-step procedure of construction of a countable family of one-particle motions for the Arratia flow. Let $(w_n)_{n\ge 0}$ be a sequence of independent standard Wiener processes, $w_k(0)=0,k\ge 0.$ Consider an arbitrary sequence $(u_n)_{n\ge 1}\colon u_n\searrow 0, n\rightarrow \infty,\, u_0\equiv 0.$ Put $\tilde{y}(u_1,t)=w_1(t)+u_1,t\in[0;1].$ Further, proceed recursively as follows. Suppose   $\tilde{y}(u_1),\ldots,\tilde{y}(u_{n-1})$ have already been built and define, for $t\in[0;1],$
\begin{equation}
\label{eq:y_tilde}
\tilde{y}(u_n,t)=\left( u_n+w_n(t) \right)\cdot \mathbbm{1}_{\{t<\tau[u_n+w_n,\tilde{y}(u_{n-1})]\}}+\tilde{y}(u_{n-1},t)\cdot \mathbbm{1}_{\{ t \ge \tau[u_n+w_n,\tilde{y}(u_{n-1})]\}}.
\end{equation}
Finally, we define
 \begin{equation*}
\left\{ {
\begin{aligned}
& \Tilde{\Tilde{y}}(u_n,t)=\tilde{y}(u_n,t)\cdot\mathbbm{1}_{\{t<\tau[\tilde{y}(u_n),w_0]\}}+w_0(t)\cdot \mathbbm{1}_{\{t\ge\tau[\tilde{y}(u_n),w_0]\}} , \  t\in[0;1],\, n\in\mathbb{N}, \\
& \Tilde{\Tilde{y}}(0,t)=w_0(t),  \  t\in[0;1].
\end{aligned} }
\right.
\end{equation*}
It follows from the definition of the Arratia flow that the sequence $(\Tilde{\Tilde{y}}(u_n))_{n\ge 0}$ has the same distribution as that of $(y(u_n))_{n\ge 0},$ since these two processes in $\mathrm{C}([0;1])^{\infty}$ have the same finite-dimensional distributions. We will say $(\tilde{\tilde{y}}(u_n))_{n\ge 0}$ is derived from $(w_n)_{n\ge 0}.$

Fix $\ve>0.$ Put $\psi\colon t\mapsto \sqrt{2t\ln\ln t^{-1}}.$
Consider  some $\alpha\in (0;1),$ whose choice will be specified later, and put $t_n=\alpha^n, \, n\geq n_0,$ where $n_0$ is such that, for all $n\geq n_0,$ it is possible to define $u_n=(1-\ve)\psi(t_n);$
 in the sequel, only such $n$ are considered. 
Let $({w}_n)_{n \ge n_0}$ be a sequence of independent standard Wiener processes that start from $0$ and are independent of $y.$ Let $(\tilde{\tilde{y}})_{n \geq n_0}$ be derived from $(w_n)_{n\ge 0}$ in compliance with the procedure described above.

Define, for $n\geq n_0,$
\begin{equation*}
\begin{aligned}
 & B_n=\left\{\frac{\nu(t_n)}{\psi(t_n)}\ge 1-\ve\right\}=\{\tau[y(u_n),y(0)]\le t_n\}, \\
 & C_n=\left\{ \tau[u_n+w_n,w_0]\le t_n \right\}, \\
 & \tilde{B}_n=\left\{ \tau[\tilde{\tilde{y}}(u_n),\tilde{\tilde{y}}(0)]\le t_n \right\}. 
 \end{aligned}
\end{equation*}
Since the sequences $(\Tilde{\Tilde{y}}(u_n))_{n\ge n_0}$ and $(y(u_n))_{n\ge n_0}$ coincide in distribution,
we have
$$
\Prob\left( \limsup_{n\rightarrow \infty} \tilde{B}_n\right)=\Prob\left( \limsup_{n\rightarrow \infty} {B}_n\right).
$$

We claim  that
\begin{equation*}
 C_n\backslash \tilde{B}_n \subset \left\{ \tau[\tilde{y}(u_{n-1}), \tilde{y}(u_n)] \le t_n\right\}, n\geq n_0+1,
\end{equation*}
where $\tilde{y}(u_{n-1})$ is defined via~\eqref{eq:y_tilde}. Indeed,
\begin{multline*}
  C_n\backslash \tilde{B}_n =\{\tau[u_n+w_n,w_0]\le t_n, \tau[\tilde{\tilde{y}}(u_n),\tilde{\tilde{y}}(0)]> t_n\}=\\
  =\{\tau[u_n+w_n,w_0]\le t_n, \tau[\tilde{\tilde{y}}(u_n),w_0]> t_n\}=\\
=\{\tau[u_n+w_n,w_0]\le t_n, \tau[\tilde{\tilde{y}}(u_n),w_0]> t_n, \tau[\tilde{y}(u_{n-1}),\tilde{y}(u_n)]\leq t_n\}.
\end{multline*}

Thus, if
\begin{multline}
 \label{eq:p(A_n)_finite}
\sum_{n\ge n_0+1} \Prob\left\{ \tau[\tilde{y}(u_{n-1}), \tilde{y}(u_n)] \le t_n \right\}=
\sum_{n\ge n_0+1} \Prob\left\{ \tau[\tilde{y}(u_{n-1}), u_n+w_n] \le t_n \right\}=\\
=\sum_{n\ge n_0+1} \Prob\left\{ \tau[u_{n-1}+w_{n-1}, u_n+w_n] \le t_n \right\}<+\infty,
\end{multline}
then
$\Prob\left( \limsup_{n\rightarrow \infty} C_n\right)=1$ implies  $\Prob\left( \limsup_{n\rightarrow \infty} {B}_n\right)=1.$

We now verify whether~\eqref{eq:p(A_n)_finite} holds. On this way, note that
$$
\Prob\left\{ \tau[u_{n-1}+w_{n-1}, u_n+w_n] \le t_n \right\}=\Prob\left\{ \theta(\frac{u_{n-1}-u_{n}}{\sqrt{2}}) \le t_n \right\},
$$
so
\begin{multline*}
\sum_{n\ge n_0+1} \Prob\left\{ \tau[\tilde{y}(u_{n-1}), y(u_n)] \le t_n \right\}= \sum_{n\geq n_0+1} \sqrt{\frac{2}{\pi}} \int_{\frac{u_{n-1}-u_{n}}{\sqrt{2t_n}}}^{+\infty}  \e^{-\frac{v^2}{2}}dv \le \\
\leq \frac{2}{\sqrt{\pi}} \sum_{n\ge n_0+1} \frac{\sqrt{t_n}}{u_{n-1}-u_n} \e^{-\frac{(u_{n-1}-u_n)^2}{4t_n}}.
\end{multline*}
Here,
\begin{multline*}
 u_{n-1}-u_n=(1-\ve)\sqrt{2\ln n} \, \alpha^{\frac{n}{2}}\left( \sqrt{\frac{\alpha^{n-1}\ln\ln \alpha^{-(n-1)}}{\alpha^{n}\ln n}}-\sqrt{\frac{\ln\ln \alpha^{-n}}{\ln n}} \right)=\\
=(1-\ve)\sqrt{2\ln n} \, \alpha^{\frac{n}{2}}\left( \sqrt{\frac{1}{\alpha}}\sqrt{\frac{\ln((n-1)\ln \alpha^{-1})}{\ln n}}-\sqrt{\frac{\ln (n \ln\alpha^{-1})}{\ln n}} \right).
\end{multline*}
As
\begin{equation*}
\label{eq:4}
\lim_{n\rightarrow\infty} \left( \sqrt{\frac{1}{\alpha}}\sqrt{\frac{\ln((n-1)\ln \alpha^{-1})}{\ln n}}-\sqrt{\frac{\ln (n \ln\alpha^{-1})}{\ln n}} \right)=\frac{1}{\sqrt{\alpha}}-1,
\end{equation*}
 the series in~\eqref{eq:p(A_n)_finite} converges if the following series converges:
$$
\sum_{n\ge n_0+1} \frac{1}{\sqrt{\ln n}} \exp\Big\{{-\frac{(1-\ve)^2 \ln n \left( \sqrt{\frac{1}{\alpha}}\sqrt{\frac{\ln((n-1)\ln \alpha^{-1})}{\ln n}}-\sqrt{\frac{\ln (n \ln\alpha^{-1})}{\ln n}} \right)^2 }{2}}\Big\},
$$
which holds whenever $\alpha$ is such that
\begin{equation}
\label{eq:alpha_condition}
\frac{(1-\ve)^2}{2} \left( \frac{1}{\sqrt{\alpha}}-1\right)^2>1.
\end{equation}

So, it is left to prove
$$
\Prob\left( \limsup_{n\rightarrow \infty} C_n\right)=1
$$
or, equivalently,
$$
\Prob\left( \liminf_{n\rightarrow \infty} {C_n^{-}}\right)=0,
$$
where the sign $\phantom{\vline}^{-}$ denotes the complement. The latter holds if
\begin{equation}
 \label{eq:lim_C_C=0}
\lim_{n\rightarrow \infty}\lim_{N\rightarrow +\infty}\Prob\left( \bigcap_{k=n}^{N} {C_k^{-}}\right)=0.
\end{equation}

Thus, we have to prove~\eqref{eq:lim_C_C=0}. For the rest of the proof, the number $n$ is fixed. Here,
\begin{multline*}
 \bigcap_{k=n}^{N} {C_k^{-}}=
\bigcap_{k=n}^{N}\{\tau[u_k+w_k,w_0]> t_k\}
=\left\{ \max_{t\in[0;t_k]}(w_0(t)-w_k(t))<u_k,\, k=\overline{n,N} \right\}=
\\
=\left\{\sup_{t\in T^{n,N}} X_t^{n,N}<1 \right\},
\end{multline*}
where $\left\{ X^{n,N}_t \mid t\in T^{n,N}\right\}$ is a centered Gaussian process:
\begin{equation*}
\left\{ {
\begin{aligned}
 T^{n,N} =& \bigcup_{k=n}^{N} T_k^{n,N}, \, T_n^{n,N} = [0;t_n],\, T_k^{n,N}=[k-1;k-1+t_k], k=\overline{n+1,N}, \\
 X_t^{n,N} =& \mathbbm{1}_{\{t\in T^{n,N}_n\}}\cdot \frac{w_0(t)-w_n(t)}{u_n}+ \\
 +& \sum_{k=n+1}^{N} \mathbbm{1}_{\{t\in T_k^{n,N}\}}\cdot \frac{w_0(t-(k-1))-w_k(t-(k-1))}{u_k}.
\end{aligned} }
\right.
\end{equation*}
Note that the sets $T_k^{n,N},k=\overline{n,N}$ are disjoint, since $t_n<1,n\geq n_0.$

Denote $\sup_{t\in T^{n,N}} X^{n,N}_t$ by $\xi^{n,N}.$ Obviously, $\xi^{n,N}< \infty$ a.s.. Further,
\begin{multline}
 \label{eq:sup_D_x}
\sigma^{n,N}=\sup_{t\in T^{n,N}} \mathrm{Var} (X^{n,N}_t)=
\max_{k=\overline{n,N}} \frac{\mathrm{Var}(w_0(t_k)-w_n(t_k))}{u_k^2}=
\max_{k=\overline{n,N}} \frac{2t_k}{u_k^2}=\\
=\max_{k=\overline{n,N}} \frac{1}{(1-\epsilon)^2 \ln\ln t_k^{-1}}=\frac{1}{(1-\epsilon)^2 \left(\ln n+ \ln\ln \alpha^{-1}\right)},
\end{multline}
so $\sigma^{n,N}\equiv\sigma^n$ does not depend on $N.$

On the set $T^{n,N},$ consider a pseudometric $\rho_{X^{n,N}}$  induced by the process $X^{n,N}:$
$$
\rho_{X^{n,N}}(t,s)=\sqrt{\mathrm{Var}(X^{n,N}_t-X^{n,N}_s)}.
$$
Estimating the metrical capacity~\cite{lifshitz,ledoux} $M^{n,N}(\delta)$ of the set $T^{n,N}$ from below,
we claim that 
 for all sufficiently large $N$
 there exists $\delta=\delta^{n,N}$ such that
\begin{equation}
\label{eq:capacity}
M^{n,N}(\delta^{n,N})\geq N-n.
\end{equation}

For that, denote $\sup_{t\in T^{n,N}_k} t$ by $s_k,\, k=\overline{n,N},$ and check that $\delta^{n,N}$ can be chosen in such a way that a set $\{ s_k \mid k=\overline{n,N}\}$ is a $\delta^{n,N}-$distinguishable subset of $T^{n,N}.$ Indeed, if $j<k,$
\begin{multline*}
 \left( \rho_{X^{n,N}}(s_j,s_k) \right)^2=\mathrm{Var} \left( \frac{w_0(s_k)-w_k(s_k)}{u_k}-\frac{w_0(s_j)-w_j(s_j)}{u_j} \right)= \\
= \mathrm{Var}\left(\frac{w_0(s_k)}{u_k}-\frac{w_0(s_j)}{u_j}\right)
+\mathrm{Var}\left( \frac{w_k(s_k)}{u_k} \right)+\mathrm{Var}\left( \frac{w_j(s_j)}{u_j} \right)= \\
= \frac{1}{u_k^2u_j^2}\mathrm{Var}\Big( u_j(w_0(s_k)-w_0(s_j))+(u_j-u_k)w_0(s_j) \Big)
+\frac{s_k}{u_k^2}+\frac{s_j}{u_j^2}=\\
=\frac{1}{u_k^2u_j^2}\left(
u_j^2 (s_k-s_j)+(u_k-u_j)^2 s_j
\right)+\frac{s_k}{u_k^2}+\frac{s_j}{u_j^2}>\frac{s_k}{u_k^2}+\frac{s_j}{u_j^2}.
\end{multline*}
Thus, it is sufficient to choose $\delta^{n,N}$ such that
\begin{multline*}
\left(\delta^{n,N}\right)^2\leq \min_{k,j=\overline{n,N},k>j} \left( \frac{s_k}{u_k^2}+\frac{s_j}{u_j^2} \right)= \\ 
= \frac{1}{2(1-\ve)^2}
\min_{k,j=\overline{n,N},k>j}\left( \frac{1}{\ln\ln\alpha^{-k}}+\frac{1}{\ln\ln\alpha^{-j}} \right).
\end{multline*}
$\delta^{n,N}=\frac{1}{1-\ve}\sqrt{\frac{1}{\ln N+\ln\ln\alpha^{-1}}}$ satisfies this condition, and, for such $\delta^{n,N}$ and for all $s,t\in \{ s_k \mid k=\overline{n,N}\}, s\neq t,$
$$
\delta^{n,N}\leq   \rho_{X^{n,N}}(s,t).
$$
This proves~\eqref{eq:capacity}.

A direct application of the Sudakov inequality~\cite{lifshitz} yields
\begin{equation}
 \label{eq:sudakov}
\E\xi^{n,N}\geq \left( 1 - \frac{1}{\sqrt{2M^{n,N}(\delta^{n,N})}} \right) \cdot \delta^{n,N} \sqrt{\ln M^{n,N}(\delta^{n,N})},
\end{equation}
 the estimate holding for $N-n\geq 24$ \cite{lifshitz}.

The concentration inequality for the Gaussian measure~\cite{ledoux}[Equation 1.23] gives that, for any $r>0,$
\begin{equation}
 \label{ledoux}
\Prob\left\{ \xi^{n,N} \leq\E \xi^{n,N} -r \right\}\le \e^{-\frac{r^2}{2\sigma^{n,N}}},
\end{equation}
where $\sigma^{n,N}=\sigma^n$ does not depend on $N$ by~\eqref{eq:sup_D_x}. 

Now fix some $r\geq 0$ and $\beta\in (0;1)$ such that
$$
\frac{1-\beta}{1-\ve}=1+r.
$$
Then, as
$$
\frac{1}{\sqrt{2M^{n,N}(\delta^{n,N})}}\leq \frac{1}{\sqrt{2(N-n)}},
$$
and
$$
\delta^{n,N} \sqrt{\ln M^{n,N}(\delta^{n,N})}\geq\frac{1}{1-\ve}\sqrt{\frac{\ln (N-n)}{\ln N+\ln\ln\alpha^{-1}}},
$$
there exists $N_0$ such that, for all $N\geq N_0,$ by\eqref{eq:sudakov},
$$
\E\xi^{n,N}> \frac{1-\beta }{1-\ve}=1+r.
$$
It follows that
$$
1< \E\xi^{n,N}-r.
$$
Recalling
$$
\bigcap_{k=n}^{N} {C_k^{-}}=\left\{ \xi^{n,N}<1 \right\},
$$
we have, by \eqref{ledoux},
$$
\Prob\left(\bigcap_{k=n}^{N} {C_k^{-}}\right)=\Prob\left\{ \xi^{n,N}<1 \right\}\leq
\Prob\left\{ \xi^{n,N}< \E\xi^{n,N}-r \right\}\leq \e^{-\frac{r^2}{2\sigma^{n}}}.
$$
Thus, as $\sigma^n\rightarrow 0,$  $n\rightarrow \infty,$
$$
\lim_{n\rightarrow \infty}\lim_{N\rightarrow +\infty}\Prob\left( \bigcap_{k=n}^{N} {C_k^{-}}\right)=0.
$$

This finishes  the proof.

\begin{rem}
In~\cite{Dorogovtsev_mono}, the Arratia flow with a drift is defined as a system of coalescing particles, each of which performs a drifted Brownian motion; the rule of coalescing is  the same: independence before the meeting and merging after. It is proved that if the drift satisfies the Lipschitz condition on the real axis and is bounded, a flow exists, and its distribution in the specified space is absolute continuous with respect to the distribution of the Arratia flow with a zero  drift (actually, the result can be extended to the case of an unbounded Lipschitz drift with some additional calculations). As an application of this analogue of the Girsanov theorem, we obtain that the statement of Theorem~\ref{theor:sizes_clusters_arratia} holds for any Arratia flow with a Lipschitz drift.
\end{rem}


\end{document}